\documentclass[11pt,reqno]{amsart}
\usepackage{amsmath,amssymb,graphicx,amscd}
\usepackage{mathrsfs}
\usepackage{enumerate}
\usepackage{mathrsfs}

\newtheorem{thm}{Theorem}
\newtheorem{cor}[thm]{Corollary}

\newtheorem{prop}[thm]{Proposition}
\theoremstyle{definition}
\newtheorem{defn}[thm]{Definition}
\theoremstyle{remark}
\newtheorem{rem}[thm]{Remark}
\numberwithin{equation}{section}
\title{Filtrations}
\author{Delia Coculescu}
\address{ETHZ \\ Departement Mathematik, R\"{a}mistrasse 101\\ Z\"{u}rich 8092, Switzerland}
 \email{delia.coculescu@math.ethz.ch}
\author{Ashkan Nikeghbali}
\address{Institut f\"ur Mathematik \\
 Universit\"at Z\"urich \\
 Winterthurerstrasse 190 CH-8057 Z\"urich \\
 Switzerland}
\email{ashkan.nikeghbali@math.unizh.ch}

\begin{document}
\begin{abstract}
In this article, we define the notion of a filtration and then give the basic theorems on initial and progressive enlargements of filtrations.
\end{abstract}

\maketitle
\section*{Definitions}
Filtrations have been introduced by Doob and have been  a fundamental feature of the theory of stochastic processes. Most basic objects, such as martingales, semimartingales, stopping times or Markov processes involve the notion of filtration.
\begin{defn}
Let $(\Omega,\mathcal{F},\mathbb{P})$ be a probability space. A \textit{filtration} on $(\Omega,\mathcal{F},\mathbb{P})$ is an increasing family $(\mathcal{F}_t)_{t\geq0}$ of sub-$\sigma$-algebras of $\mathcal{F}$. In other words, for each $t$,  $\mathcal{F}_t$ is a $\sigma$-algebra included in $\mathcal{F}$ and if $s\leq t$, $\mathcal{F}_s\subset\mathcal{F}_t$. A probability space $(\Omega,\mathcal{F},\mathbb{P})$ endowed with a filtration $(\mathcal{F}_t)_{t\geq0}$ is called a filtered probability space.
\end{defn}
We now give a definition which is very closely related to that of a filtration:
\begin{defn}
A stochastic process $(X_t)_{t\geq0}$ on $(\Omega,\mathcal{F},\mathbb{P})$ is adapted to the filtration $(\mathcal{F}_t)$ if, for each $t\geq0$, $X_t$ is $\mathcal{F}_t$-measurable. 
\end{defn}
A stochastic process $X$ is always adapted to its \textit{natural filtration} $\mathcal{F}_t^X=\sigma(X_s,\;s\leq t)$ (the last notation meaning that $\mathcal{F}_t$ is the smallest $\sigma$-algebra with respect to which all the variables $(X_s,\;s\leq t)$ are measurable). $(\mathcal{F}_t^X)$ is hence the smallest filtration to which $X$ is adapted.

The parameter $t$ is often thought of  as time, and the $\sigma$-algebra $\mathcal{F}_t$ represents the set of information available at time $t$, that is events that have occurred up to time $t$. The filtration $(\mathcal{F}_t)_{t\geq0}$ thus represents the evolution of the information or knowledge of the world with time. If $X$ is an adapted process, then $X_t$, its value   at time $t$, only depends on the evolution of the universe prior to $t$.

\begin{defn}
Let $\left( \Omega ,\mathcal{F},\left( \mathcal{F}_{t}\right) _{t\geq 0},%
\mathbb{P}\right) $ be a filtered probability space. 
\begin{enumerate}[(i)]
\item The filtration $\left( \mathcal{F}_{t}\right) _{t\geq 0}$ is said to be\textit{ complete} if $\left( \Omega ,\mathcal{F},
\mathbb{P}\right) $ is complete and if $\mathcal{F}_0$ contains all the $\mathbb{P}$-null sets.
\item The filtration $\left( \mathcal{F}_{t}\right) _{t\geq 0}$ is said to satisfy the \textit{usual hypotheses} if it is complete and right continuous, that is $\mathcal{F}_t=\mathcal{F}_{t+}$, where $$\mathcal{F}_{t+}=\bigcap_{u>t}\mathcal{F}_u.$$
\end{enumerate}
Some fundamental theorems, such as the Debut theorem, require the usual hypotheses. Hence naturally, very often in the literature on the theory of stochastic processes and mathematical finance, the underlying filtered probability spaces are assumed to satisfy the usual hypotheses. This assumption is not very restrictive for the following reasons:

\begin{enumerate}[(a)]
\item Any filtration can  easily be made complete and right continuous; indeed, given a filtered probability space $\left( \Omega ,\mathcal{F},\left( \mathcal{F}_{t}\right) _{t\geq 0},%
\mathbb{P}\right) $, we first complete the probability space $\left( \Omega ,\mathcal{F},\mathbb{P}\right) $, and then we add all the $\mathbb{P}$-null sets to every $\mathcal{F}_{t+}$, $t\geq0$. The new filtration thus obtained  satisfies the usual hypotheses and is called the usual augmentation of $\left( \mathcal{F}_{t}\right) _{t\geq 0}$;
\item Moreover, in most classical and encountered cases, the filtration  $\left( \mathcal{F}_{t}\right) _{t\geq 0}$ is right continuous. Indeed, this is the case when for instance  $\left( \mathcal{F}_{t}\right) _{t\geq 0}$ is the natural filtration of a Brownian Motion, a L\'evy process, a Feller process or a Hunt process (see \cite{protter,revuzyor}).
\end{enumerate}

\end{defn}
\section*{Enlargements of filtrations}
For more precise and detailed references, the reader can consult the books \cite{jeulin,jeulinyor,columbia,protter} or the survey article \cite{ashkanessay}.
\subsection*{Generalities}
Let $\left( \Omega ,\mathcal{F},\left( \mathcal{F}_{t}\right) _{t\geq 0},%
\mathbb{P}\right) $ be a filtered probability space satisfying the usual hypotheses. Let $\left( \mathcal{G}_{t}\right) _{t\geq 0}$ be another filtration satisfying the usual hypotheses and such that $\mathcal{F}_{t}\subset\mathcal{G}_{t}$ for every $t\geq0$. One natural question is: how are the $\left( \mathcal{F}_{t}\right)$-semimartingales modified when considered as stochastic processes in the larger filtration $\left( \mathcal{G}_{t}\right)$? Given the importance of semimartingales and martingales (in particular in mathematical finance where they are used to model prices), it seems natural to characterize  situations where the semimartingale or martingale properties are preserved:
\begin{defn}
We shall say that the pair of filtrations $\left(\mathcal{F}_{t},
\mathcal{G}_{t}\right)$ satisfies the $\left(H'\right)$ hypothesis
if every $\left(\mathcal{F}_{t}\right)$-semimartingale is a
$\left(\mathcal{G}_{t}\right)$-semimartingale.
\end{defn}
\begin{rem}
In fact, using a classical decomposition of semimartingales due to Jacod and M\'emin, it is enough to check that every
$\left(\mathcal{F}_{t}\right)$-bounded martingale is a
$\left(\mathcal{G}_{t}\right)$-semimartingale.
\end{rem}
\begin{defn}
We shall say that the pair of filtrations $\left(\mathcal{F}_{t},
\mathcal{G}_{t}\right)$ satisfies the $\left(H\right)$ hypothesis
if every $\left(\mathcal{F}_{t}\right)$-local martingale is a
$\left(\mathcal{G}_{t}\right)$-local martingale.
\end{defn}

 The techniques to answer such questions have been developed in the late 70's under the name of the theory of enlargements of filtrations. The theory of enlargements of filtrations has been recently very widely used in mathematical finance, specially in  insider trading models and even more spectacularly in models of default risk.  The insider trading models are usually based on the so called \textit{initial enlargements  of filtrations} whereas the models of default risk fit perfectly well in the framework of the \textit{progressive enlargements of filtrations}.  More precisely, given a filtered probability space
$\left(\Omega,\mathcal{F},\left(\mathcal{F}_{t}\right),\mathbb{P}\right)$,
there are essentially two ways of enlarging filtrations:
\begin{itemize}
\item \textit{initial enlargements}, for which
$\mathcal{G}_{t}=\mathcal{F}_{t}\bigvee\mathcal{H}$, i.e. the new
information $\mathcal{H}$ is brought in at the origin of time; and
\item \textit{progressive enlargements}, for which
$\mathcal{G}_{t}=\mathcal{F}_{t}\bigvee\mathcal{H}_{t}$, i.e. the
new information is brought in progressively as the time $t$
increases.
\end{itemize}
Before presenting the basic theorems on enlargements of filtrations, we state a useful theorem due to Stricker:

\begin{thm}[Stricker \cite{stricker}]
Let $\left(\mathcal{F}_{t}\right)$ and
$\left(\mathcal{G}_{t}\right)$ be two filtrations as above, such that for all
$t\geq0$, $\mathcal{F}_{t}\subset\mathcal{G}_{t}$. If
$\left(X_{t}\right)$ is a $\left(\mathcal{G}_{t}\right)$
semimartingale which is $\left(\mathcal{F}_{t}\right)$ adapted, then
it is also an $\left(\mathcal{F}_{t}\right)$ semimartingale.
\end{thm}
\subsection*{Initial enlargements of filtrations}
The most important theorem on initial enlargements of filtrations is due to Jacod and deals with the special case where the initial information brought in at the origin of time consists of the $\sigma$-algebra generated by a random variable. More precisely
let
$\left(\Omega,\mathcal{F},\left(\mathcal{F}_{t}\right),\mathbb{P}\right)$
be a filtered probability space satisfying the usual assumptions.
Let $Z$ be an $\mathcal{F}$ measurable random variable. Define
$$\mathcal{G}_{t}=\bigcap_{\varepsilon>0}\left(\mathcal{F}_{t+\varepsilon}\bigvee\sigma\left\{Z\right\}\right).$$In financial models, the filtration $\left(\mathcal{F}_{t}\right)$ represents the public information in a financial market and the random variable $Z$ stands for the additional (anticipating) information of an insider.

The conditional laws of $Z$ given $\mathcal{F}_{t}$, for $t\geq0$
play a crucial role in initial enlargements.
\begin{thm}[Jacod's criterion]\label{thmdejacod}
Let $Z$ be an $\mathcal{F}$ measurable random variable and let
$Q_{t}\left(\omega,dx\right)$ denote the regular conditional
distribution of $Z$ given $\mathcal{F}_{t},\;t\geq0$. Suppose that
for each $t\geq0$, there exists a positive $\sigma$-finite measure
$\eta_{t}\left(dx\right)$ (on
$\left(\mathbb{R},\mathcal{B}\left(\mathbb{R}\right)\right)$) such
that
$$Q_{t}\left(\omega,dx\right)\ll\eta_{t}\left(dx\right)\;\mathrm{a.s.}$$Then
every $\left(\mathcal{F}_{t}\right)$-semimartingale is a
$\left(\mathcal{G}_{t}\right)$-semimartingale.
\end{thm}
\begin{rem}
In fact this theorem still holds for random variables with values in
a standard Borel space. Moreover, the existence of the
$\sigma$-finite measure $\eta_{t}\left(dx\right)$ is equivalent to
the existence of one positive $\sigma$-finite measure
$\eta\left(dx\right)$ such that
$Q_{t}\left(\omega,dx\right)\ll\eta\left(dx\right)$ and in this case
$\eta$ can be taken to be the distribution of $Z$.
\end{rem}
Now we give classical corollaries of Jacod's theorem.
\begin{cor}
Let $Z$ be independent of $\mathcal{F}_{\infty}$. Then every
$\left(\mathcal{F}_{t}\right)$-semimartingale is a
$\left(\mathcal{G}_{t}\right)$-semimartingale.
\end{cor}
\begin{cor}
Let $Z$ be a random variable taking on only a countable number of
values. Then every $\left(\mathcal{F}_{t}\right)$-semimartingale is
a $\left(\mathcal{G}_{t}\right)$-semimartingale.
\end{cor}
It is possible to obtain in some cases an explicit decomposition of an $\left(\mathcal{F}_{t}\right)$-local martingale as a $\left(\mathcal{G}_{t}\right)$-semimartingale (see \cite{jeulin, jeulinyor, columbia, ashkanessay, protter}). For example, if $Z=B_{t_0}$, for some fixed time $t_0>0$ and a  Brownian Motion $B$, it can be shown that Jacod's criterion holds for $t<t_0$ and that every $\left(\mathcal{F}_{t}\right)$-local martingale is a semimartingale for $0\leq t<t_0$, but not necessarily including $t_0$. There are indeed in this case $\left(\mathcal{F}_{t}\right)$-local martingales which are not $\left(\mathcal{G}_{t}\right)$-semimartingales. Moreover, $B$ is a $\left(\mathcal{G}_{t}\right)$-semimartingale which decomposes as:
$$B_{t}=B_{0}+\widetilde{B}_{t}+\int_{0}^{t\wedge t_{0}}ds\dfrac{B_{t_{0}}-B_{s}}{t_{0}-s},$$where $\left(\widetilde{B}_{t}\right)$ is a $\left(
\mathcal{G}_{t}\right)$ Brownian Motion. 
\begin{rem}
There are important cases where Jacod's criterion does not hold but where other methods apply (\cite{jeulin, columbia, ashkanessay}.
\end{rem}
\subsection*{Progressive enlargements of filtrations}
Let $\left( \Omega ,\mathcal{F},\left( \mathcal{F}_{t}\right) _{t\geq 0},%
\mathbb{P}\right) $ be a filtered probability space satisfying the usual hypotheses, and $\rho :$
$\left(
\Omega ,\mathcal{F}\right) \rightarrow \left( \mathbb{R}_{+},\mathcal{B}%
\left( \mathbb{R}_{+}\right) \right) $ be a random time. We enlarge
the initial filtration $\left( \mathcal{F}_{t}\right) $\ with the
process $\left( \rho \wedge t\right) _{t\geq 0}$, so that the new
enlarged filtration $\left( \mathcal{F}_{t}^{\rho }\right) _{t\geq
0}$\ is the
smallest filtration (satisfying the usual assumptions) containing $\left( \mathcal{F}_{t}\right) $\ and making $%
\rho $\ a stopping time (i.e. $\mathcal{F}_{t}^{\rho
}=\mathcal{K}_{t+}^{o}$, where
$\mathcal{K}_{t}^{o}=\mathcal{F}_{t}\bigvee\sigma\left(\rho\wedge
t\right)$). One may interpret  $\rho$ as the instant of default of an issuer; the 
given filtration $(\mathcal{F}_t)$ can be thought of as the filtration of
default-free prices, for which $\rho$ is not a
stopping time. Then, the filtration
$(\mathcal{F}_t^\rho)$ is the defaultable market filtration used for the pricing of defaultable assets.

A few processes will play a crucial role in
our discussion:

\begin{itemize}
\item the $\left( \mathcal{F}_{t}\right) $-supermartingale
\begin{equation}
Z_{t}^{\rho }=\mathbb{P}\left[ \rho >t\mid \mathcal{F}_{t}\right]
\label{surmart}
\end{equation}%
chosen to be c\`{a}dl\`{a}g, associated to $\rho $\ by Az\'{e}ma;

\item the $\left( \mathcal{F}_{t}\right) $-dual optional
projection of the process $1_{\left\{ \rho \leq t\right\} }$,
denoted by $A_{t}^{\rho }$;

\item the c\`{a}dl\`{a}g martingale
\begin{equation*}
\mu _{t}^{\rho }=\mathbb{E}\left[ A_{\infty }^{\rho }\mid \mathcal{F}_{t}%
\right] =A_{t}^{\rho }+Z_{t}^{\rho }.
\end{equation*}%
\end{itemize}
\begin{thm}
Every $\left( \mathcal{F}_{t}\right) $-local
martingale $\left( M_{t}\right) $, stopped at $\rho $, is a $\left( \mathcal{%
F}_{t}^{\rho }\right) $\-semimartingale, with canonical
decomposition:
\begin{equation}
M_{t\wedge \rho }=\widetilde{M}_{t}+\int_{0}^{t\wedge \rho
}\dfrac{d\langle M,\mu ^{\rho }\rangle_{s}}{Z_{s-}^{\rho }}
\label{decocanonique}
\end{equation}%
where $\left( \widetilde{M}_{t}\right) $\ is an $\left( \mathcal{F}%
_{t}^{\rho }\right) $-local martingale.
\end{thm}
The most interesting case in the theory of progressive enlargements
of filtrations is when $\rho$ is an honest time or equivalently the end of an $\left(
\mathcal{F}_{t}\right) $\ optional set $\Gamma $, i.e\textbf{\ }%
\begin{equation*}
\rho=\sup \left\{ t:\left( t,\omega \right) \in \Gamma \right\}.
\end{equation*}%
Indeed, in this case, the pair of filtrations $(\mathcal{F}_t,\mathcal{F}_t^\rho)$ satisfies the $(H')$ hypothesis: every $\left( \mathcal{F}_{t}\right) $-local
martingale $\left( M_{t}\right) $, is an $\left( \mathcal{%
F}_{t}^{\rho }\right) $-semimartingale, with canonical
decomposition:
\begin{equation*}
M_{t}=\widetilde{M}_{t}+\int_{0}^{t\wedge \rho }\frac{d\langle M,\mu
^{\rho}\rangle_{s}}{Z_{s-}^{\rho}}-\int_{\rho}^{t}\dfrac{d\langle M,\mu ^{\rho
}\rangle_{s}}{1-Z_{s-}^{\rho}}. 
\end{equation*}
The next decomposition formulae are widely used for pricing in default models:
\begin{prop}
\begin{enumerate}[(i)]
\item Let $\xi\in L^{1}$. Then a c\`{a}dl\`{a}g version of the
martingale
$\xi_{t}=\mathbb{E}\left[\xi|\mathcal{F}_{t}^{\rho}\right]$ is given
by:
$$\xi_{t}=\dfrac{1}{Z_{t}^{\rho}}\mathbf{1}_{t<\rho}\mathbb{E}\left[\xi\mathbf{1}_{t<\rho}|\mathcal{F}_{t}\right]+\xi\mathbf{1}_{t\geq\rho}.$$
\item Let $\xi\in L^{1}$ and let $\rho$ be an honest time. Then a c\`{a}dl\`{a}g version of the
martingale $\xi_{t}=\mathbb{E}\left[\xi|\mathcal{F}_{t}^{\rho}\right]$
is given by:
$$\xi_{t}=\dfrac{1}{Z_{t}^{\rho}}\mathbb{E}\left[\xi\mathbf{1}_{t<\rho}|\mathcal{F}_{t}\right]\mathbf{1}_{t<\rho}+\dfrac{1}{1-Z_{t}^{\rho}}\mathbb{E}\left[\xi\mathbf{1}_{t\geq \rho}|\mathcal{F}_{t}\right]\mathbf{1}_{t\geq \rho}.$$
\end{enumerate}
\end{prop}
\subsection*{The $(H)$ hypothesis}
 The $(H)$ hypothesis is sometimes presented as a no-abitrage condition in default models. Let
$\left(\Omega,\mathcal{F},\mathbb{P}\right)$ be a probability space
satisfying the usual assumptions. Let $\left(\mathcal{F}_{t}\right)$
and $\left(\mathcal{G}_{t}\right)$ be two sub-filtrations of
$\mathcal{F}$, with
$$\mathcal{F}_{t}\subset\mathcal{G}_{t}.$$Br\'emaud and Yor \cite{bremaudyor} have proven the following characterization of the $(H)$ hypothesis:
\begin{thm}
The following are
equivalent:
\begin{enumerate}
\item Every $\left(\mathcal{F}_{t}\right)$ martingale is a
$\left(\mathcal{G}_{t}\right)$ martingale;
\item For all $t\geq0$, the sigma fields $\mathcal{G}_{t}$ and
$\mathcal{F}_{\infty}$ are independent conditionally on
$\mathcal{F}_{t}$.
\end{enumerate}
\end{thm}
\begin{rem}
We shall also say that $\left(\mathcal{F}_{t}\right)$ is \textit{immersed} in
$\left(\mathcal{G}_{t}\right)$.
\end{rem} In the framework of the progressive enlargement of some filtration
$\left(\mathcal{F}_{t}\right)$ with a random time $\rho $, the $(H)$ hypothesis  is
equivalent to one of the following hypothesis:

\begin{enumerate}[(i)]
\item $\forall t$, the $\sigma $-algebras $\mathcal{F}_{\infty }$\ and $%
\mathcal{F}_{t}^{\rho }$\ are conditionally independent given $\mathcal{F}%
_{t}$.

\item For all bounded $\mathcal{F}_{\infty }$ measurable random variables $%
\mathbf{F}$\ and all bounded $\mathcal{F}_{t}^{\rho }$ measurable
random
variables $\mathbf{G}_{t}$, we have%
\begin{equation*}
\mathbb{E}\left[ \mathbf{FG}_{t}\mid \mathcal{F}_{t}\right]
=\mathbb{E}\left[ \mathbf{F}\mid \mathcal{F}_{t}\right]
\mathbb{E}\left[ \mathbf{G}_{t}\mid \mathcal{F}_{t}\right] .
\end{equation*}

\item For all bounded $\mathcal{F}_{t}^{\rho }$ measurable random variables $%
\mathbf{G}_{t}$:%
\begin{equation*}
\mathbb{E}\left[ \mathbf{G}_{t}\mid \mathcal{F}_{\infty }\right] =\mathbb{E}%
\left[ \mathbf{G}_{t}\mid \mathcal{F}_{t}\right] .
\end{equation*}

\item For all bounded $\mathcal{F}_{\infty }$ measurable random variables $%
\mathbf{F}$,
\begin{equation*}
\mathbb{E}\left[ \mathbf{F}\mid \mathcal{F}_{t}^{\rho }\right] =\mathbb{E}%
\left[ \mathbf{F}\mid \mathcal{F}_{t}\right] .
\end{equation*}

\item For all $s\leq t$,
\begin{equation*}
\mathbb{P}\left[ \rho \leq s\mid \mathcal{F}_{t}\right]
=\mathbb{P}\left[ \rho \leq s\mid \mathcal{F}_{\infty }\right] .
\end{equation*}
\end{enumerate}
Now, a natural question, specially in view of applications to financial mathematics, is: how is the $(H)$ hypothesis affected when we make an equivalent change of probability
measure?
\begin{prop}
Let $\mathbb{Q}$ be a probability measure which is equivalent to
$\mathbb{P}$ (on $\mathcal{F}$). Then every
$\left(\mathcal{F}_{\bullet},\mathbb{Q}\right)$-semimartingale is a
$\left(\mathcal{G}_{\bullet},\mathbb{Q}\right)$-semimartingale.
\end{prop}
Now, define:
$$\dfrac{d\mathbb{Q}}{d\mathbb{P}}\Big|_{\mathcal{F}_{t}}=R_{t};\quad\dfrac{d\mathbb{Q}}{d\mathbb{P}}\Big|_{\mathcal{G}_{t}}=R'_{t}.$$
If
$Y=\dfrac{d\mathbb{Q}}{d\mathbb{P}}$, then the hypothesis $(H)$
holds under $\mathbb{Q}$ if and only if:
$$\forall X\geq0,\;X\in\mathcal{F}_{\infty},\quad \dfrac{\mathbb{E}_{\mathbf{P}}\left[XY|\mathcal{G}_{t}\right]}{R'_{t}}=\dfrac{\mathbb{E}_{\mathbf{P}}\left[XY|\mathcal{F}_{t}\right]}{R_{t}}.$$
In particular, when $\dfrac{d\mathbb{Q}}{d\mathbb{P}}$ is
$\mathcal{F}_{\infty}$ measurable, $R_{t}=R'_{t}$ and the hypothesis
$(H)$ holds under $\mathbb{Q}$.

Now let us give a decomposition formula:
\begin{thm}
If $\left(X_{t}\right)$ is a
$\left(\mathcal{F}_{\bullet},\mathbb{Q}\right)$-local martingale,
then the stochastic process:
$$I_{X}\left(t\right)=X_{t}+\int_{0}^{t}\dfrac{R'_{s-}}{R'_{s}}\left(\dfrac{1}{R_{s-}}d[X,R]_{s}-\dfrac{1}{R'_{s-}}d[X,R']_{s}\right)$$is
a $\left(\mathcal{G}_{\bullet},\mathbb{Q}\right)$-local martingale.
\end{thm}

\end{document}